\documentclass[letterpaper, 10 pt, conference]{ieeeconf}  

\IEEEoverridecommandlockouts                              
\usepackage{graphics} 
\usepackage{mathtools} 
\usepackage{times} 
\usepackage{amsmath} 
\usepackage{amssymb}  
\DeclareMathOperator{\sign}{sgn}
\title{\LARGE \bf
A Deterministic Annealing Approach to the Multiple Traveling Salesmen and Related Problems
}

\author{Mayank Baranwal$^{1,a}$, Brian Roehl and Srinivasa M. Salapaka$^{1,b}$
\thanks{$^{1}$Department of Mechanical Science and Engineering, University of Illinois at Urbana-Champaign, 61801 IL, USA}
\thanks{$^{a}${\tt\small baranwa2@illinois.edu}, $^{b}${\tt\small salapaka@illinois.edu}}
}

\begin{document}

\maketitle
\thispagestyle{empty}
\pagestyle{empty}

\begin{abstract}
This paper presents a novel and efficient heuristic framework for approximating the solutions to the multiple traveling salesmen problem (m-TSP) and other variants on the TSP. The approach adopted in this paper is an extension of the Maximum-Entropy-Principle (MEP) and the Deterministic Annealing (DA) algorithm. The framework is presented as a general tool that can be suitably adapted to a number of variants on the basic TSP. Additionally, unlike most other heuristics for the TSP, the framework presented in this paper is independent of the edges defined between any two pairs of nodes. This makes the algorithm particularly suited for variants such as the close-enough traveling salesman problem (CETSP) which are challenging due to added computational complexity. The examples presented in this paper illustrate the effectiveness of this new framework for use in TSP and many variants thereof.
\end{abstract}

\section{INTRODUCTION}\label{sec:Intro}
The Traveling Salesman Problem (TSP) \cite{little1963algorithm, lin1973effective, wang2003traveling} is one of the most extensively studied optimization problems. A TSP is defined by a set of nodes or cities, and the edges connecting them which define the cost of travel between each city. Each solution to a TSP is referred as a {\em tour} and is made up of a combination of edges such that each node is visited sequentially. The optimal tour is the combination of edges that minimizes the total cost to the salesman. There are several applications of the TSP to real world problems. Common applications of the TSP are vehicle delivery route planning and toolhead path planning for VLSI circuit boards. Junger {\em et al.} \cite{junger1995traveling} as well as Bektas \cite{bektas2006multiple} have explored many more applications of the TSP to more specialized problems, demonstrating the real world value and importance of developing effective solutions to the TSP.  

The TSP belongs to the class of NP-complete problems that are computationally intensive to solve \cite{bektas2006multiple}. Generally, finding the optimal solution requires a calculation time given by $\Theta(n!)$, where $n$ is the number of nodes in the system. For relatively small data sets linear solving programs can find the optimal solution within a short period of time, but for larger data sets the computations can become extremely time intensive. Only with the development of more powerful computers has the optimal solution been discovered for larger data sets, as computation times in the hundreds of CPU-years had made those solutions infeasible in the past \cite{applegate2009certification}. Many heuristics have been developed for the TSP \cite{junger1995traveling, bektas2006multiple, johnson2007experimental} and they can offer significant runtime improvements conventional optimization for large data sets, at the sacrifice of some deviation from the optimal solution.

Despite the success of heuristics for the basic TSP, there are a quite a few variants to the TSP that can pose a challenge to some of these methodologies. Variants to the basic TSP are necessary to appropriately model realistic situations. In this paper, we present a novel framework for efficient solutions of two such variants of the TSP - a) the {\em multi-Traveling Salesmen Problem} (mTSP) allows more than one `salesman' to operate between the cities, such that the solution to the mTSP is comprised of several routes, one for each salesman, and the optimal tour would be the set of routes such that the total distance travelled is minimized. Additional constraints may be imposed on the system such as requiring each salesman to start at the same point, representing a warehouse or depot. In order to use the conventional heuristics on a mTSP, transformations can be applied in order to generate a TSP which represents the mTSP \cite{bellmore1974transformation}. b) the {\em Close Enough Traveling Salesman Problem} (CETSP) is a variant where the salesman must only come within a certain radius of each city on the tour. This adds great complexity to the problem. Because of the significant increase in the number of edges, many conventional heuristics are unable to address this variant. Special formulations have been developed to address this problem \cite{gulczynski2006close, mennell2009heuristics}.

In this paper we present a solution methodology based on Deterministic Annealing (DA) for the above variants on the TSP. DA is well-suited to combinatorial clustering/resource allocation problems that require obtaining an optimal partition of an underlying domain, and optimally assigning resources to each cell of the partition. DA-based methods have been reported in a vast number of applications such as minimum distortion problems in data compression \cite{rose1998deterministic}, model aggregation \cite{xu2014aggregation}, routing problems in multiagent networks \cite{kale2012maximum}, locational optimization problems \cite{salapaka2003constraints, salapaka2003locational}, and coverage control problems \cite{xu2014clustering}. The DA algorithm is developed at the intersection of statistical mechanics and information theory, and is a methodology that allows the greatest amount of information about a system to be inferred from on a limited amount of given information using the Maximum-Entropy-Principle (MEP). In this paper the MEP is used to consider every potential tour of the cities, and through the optimization process of DA the shortest tour through every city is determined. While the original DA algorithm was developed in the context of clustering, it was later adapted to the basic TSP as a case of constrained clustering \cite{rose1990deterministic, rose1998deterministic}, which serves as the foundation for this extension to the mTSP.

The rest of the paper is organized as follows. Section \ref{sec:Prob} presents the mathematical formulation of the mTSP and CETSP. We then provide a brief overview of the DA algorithm in Section \ref{sec:MEP}, which is then followed by the extension of DA to variants on the TSP in Section \ref{sec:methodology}. We then present a few examples to illustrate the effectiveness of this new framework for use in TSP and many variants thereof in Section \ref{sec:results}. We then conclude the paper and provide directions for future work.

\section{PROBLEM FORMULATION}\label{sec:Prob}
In this section, we describe some variants on the TSP. The set of node locations is denoted by $\{x_i:x_i\in\mathbb{R}^2, 1\leq i\leq n\}$, where $n$ is the number of nodes. We use $\alpha\in\mathbb{R}^2$ to denote the location of depot. The distance between any two nodes $i$ and $j$ is denoted by $d_{i,j}$. For example, for squared-euclidean distance, we have $d_{ij}=\|x_i-x_j\|_2^2$. We use $S_n$ to denote the set of all permutations of $\{1,\dots ,n\}$. An element of $S_n$ is denoted by $\mu$. An element $\sigma=[\sigma_1,\dots ,\sigma_n,\sigma_{n\!+\!1}]$ in $\tilde{S}_n:=S_n\cup\{n\!+\!1\}$ is given by
\begin{small}
	\begin{eqnarray}\label{eq:sigma}
		\sigma_i := \sigma(i) = \left\{
				\begin{array}{lll}
					\mu(i) & \mbox{if } &\quad i\in \{1,2,\dots ,n\} \\
					\mu(1) & \mbox{if} &\quad i = n+1.
				\end{array}
			\right.\nonumber
	\end{eqnarray}
\end{small}
\begin{figure*}[!t]
	\begin{center}
	\begin{tabular}{cccc}
		\includegraphics[width=1.5in]{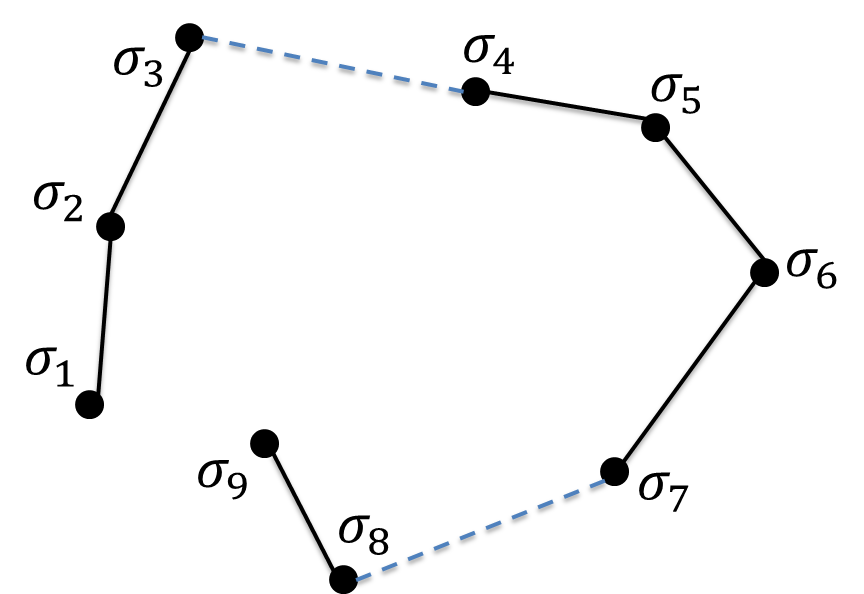}&
		\includegraphics[width=1.5in]{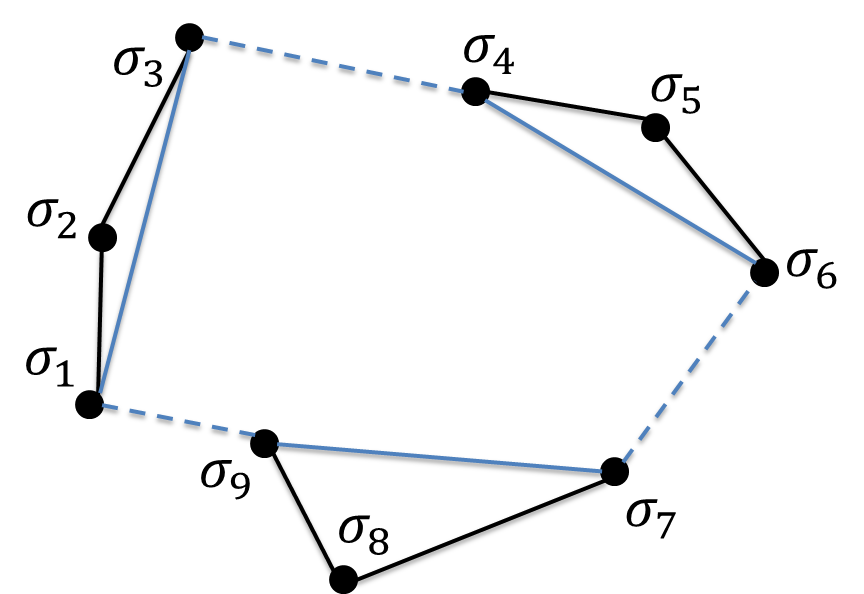}&
		\includegraphics[width=1.5in]{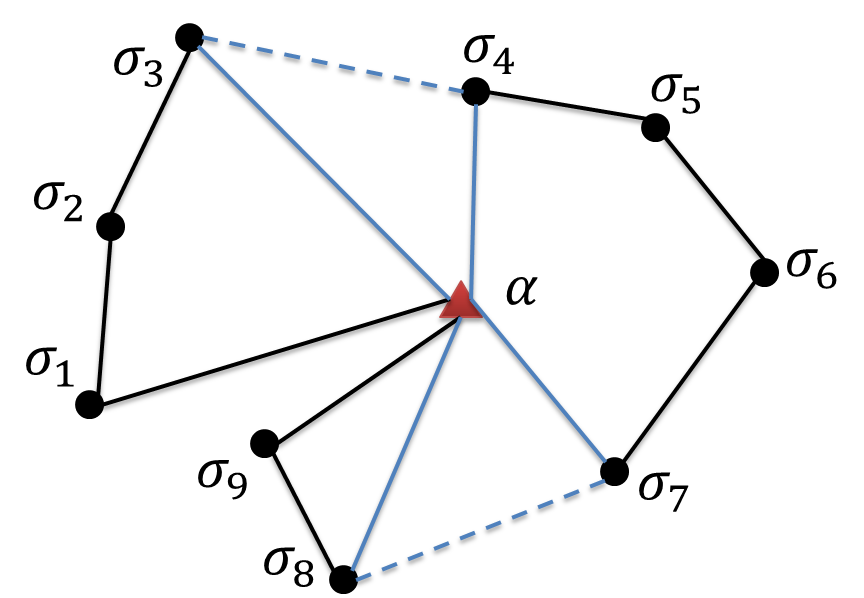}&	
		\includegraphics[width=1.5in]{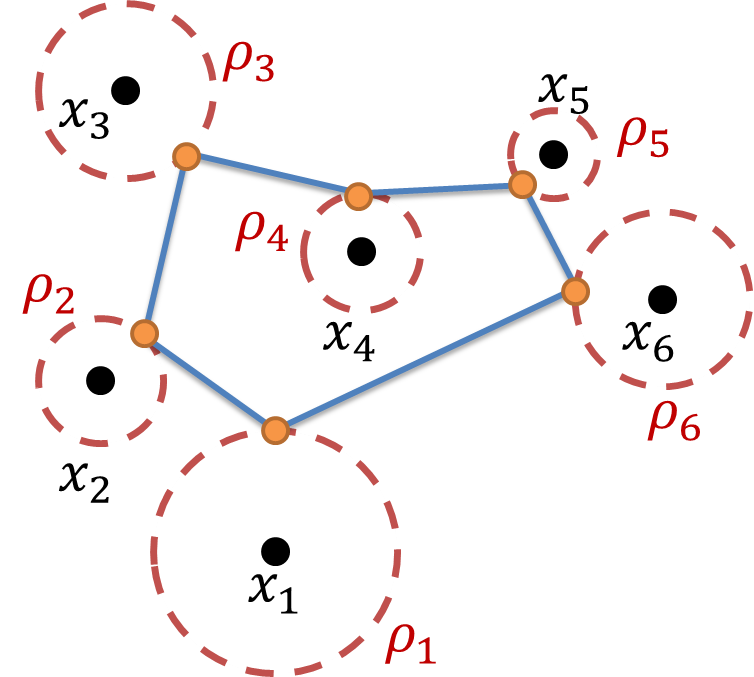}\cr
		(a)&(b)&(c)&(d)
	\end{tabular}
	\vspace{-0.5em}
	\caption{{\small Schematic of a (a) $9$-nodes non-returning $3$TSP. The dashed blue lines indicate the removal of links, which corresponds to $\alpha_{k_1=3}=1, \alpha_{k_2=7}=1$ (b) $9$-nodes returning $3$TSP. The dashed blue lines indicate the removal of links ($-d_{\sigma_{k_j},\sigma_{k_j\!+\!1}}$), while solid blue lines indicate the addition of links ($d_{\sigma_{k_j},\sigma_{k_{(j\!-\!1)}\!+\!1}}$). (c) $9$-nodes single-depot returning $3$TSP. (d) Single salesman returning CETSP. Each node $x_i$ is provided with a radius parameter $\rho_i$. The orange dots indicate $y_j$ such that $v_{ij}=1$ for some node $i$.}}
	\label{fig:mTSP_prob}
	\vspace{-2em}
	\end{center}
\end{figure*}

\subsection{Non-Returning Multi-Traveling Salesmen Problem}\label{subsec:Problem_mTSP_nonR}
In this problem, we are given a set of $n$ nodes ($\{x_i\}$) and $m$ salesmen to traverse these nodes. The objective is to minimize the total tour-length such that each node is visited just once by only one salesman. The starting and ending node locations of each salesman can not coincide. This formulation is applicable to problems pertaining to non-recurring events, such as the scheduling of orders at a steel rolling company \cite{tang2000multiple}. The non-returning $m$TSP can be mathematically described as
\begin{small}
	\begin{eqnarray}\label{eq:opt1}
		\min_{\substack{\sigma\in \tilde{S}_n\\ \{\alpha_{k_j}\}}}\left\{\sum_{i=1}^{n-1}d_{\sigma_i,\sigma_{i+1}}-\sum_{j=1}^{m-1}\alpha_{k_j}d_{\sigma_{k_j},\sigma_{k_j+1}}\right\}\nonumber \\
		\text{s.t.}\quad\forall j,\sum_{k_j}\alpha_{k_j}=1,\quad\alpha_{k_j}\in\{0,1\},\quad 1\leq k_j\leq n\!-\!1.\nonumber
	\end{eqnarray}
\end{small}
Fig. \ref{fig:mTSP_prob}a shows a schematic of a non-returning $m$TSP.

\subsection{Returning Multi-Traveling Salesmen Problem}\label{subsec:Problem_mTSP_R}
In this problem, we are given a set of $n$ nodes ($\{x_i\}$) and $m$ salesmen, and the objective is to minimize the total tour length such that each node is visited by only one salesman, and the start and end positions of each salesman must be coincident. Many recurring events, such as job scheduling \cite{gorenstein1970printing} fall under this category. The returning $m$TSP can be mathematically described as
\begin{small}
	\begin{eqnarray}\label{eq:opt2}
		\min_{\substack{\sigma\in \tilde{S}_n\\ \{\alpha_{k_j}\}}}\left\{\sum_{i=1}^{n}d_{\sigma_i,\sigma_{i+1}}\!+\!\sum_{j=1}^{m}\alpha_{k_j}\left[-d_{\sigma_{k_j},\sigma_{k_j\!+\!1}}\!+\!d_{\sigma_{k_j},\sigma_{k_{(j\!-\!1)}\!+\!1}}\right]\right\} \nonumber \\
		\text{s.t.}\quad\forall j,\sum_{k_j}\alpha_{k_j}=1,\quad\alpha_{k_j}\in\{0,1\},\quad 1\leq k_j\leq n; k_0 = k_m.\nonumber
	\end{eqnarray}
\end{small}
Fig. \ref{fig:mTSP_prob}b shows a schematic of a returning $m$TSP.

\subsection{Single-Depot Returning Multi-Traveling Salesmen Problem}\label{subsec:Problem_mTSP_RD}
In this problem, we are given a set of $n$ nodes ($\{x_i\}$), a depot ($\alpha$) and $m$ salesmen, and the objective is to determine the optimal tour such that each node is visited by only one salesman. Each salesman must start and end at the depot. The total distance traveled by all salesmen is minimized. Real-world problems such as vehicle routing problem (VRP) with single-depot fall under this category. If $x_0=\alpha$ denotes the location of depot and $\mathcal{I}$ is the indexed-set of nodes, then the single-depot $m$TSP is mathematically described as
\begin{small}
	\begin{eqnarray}\label{eq:opt1}
		\min_{\sigma\in \tilde{S}_n, \{\alpha_{k_j}\}}\bigg\{\sum_{i=1}^{n-1}d_{\sigma_i,\sigma_{i+1}}+d_{\sigma_1,\alpha}+d_{\sigma_n,\alpha} \nonumber \\
		+\!\sum_{j=1}^{m\!-\!1}\alpha_{k_j}\left[-d_{\sigma_{k_j},\sigma_{k_j+1}}+d_{\sigma_{k_j},\alpha}+d_{\sigma_{k_j+1},\alpha}\right]\bigg\} \nonumber \\
		\text{s.t.}\quad\forall j,\sum_{k_j}\alpha_{k_j}=1,\quad\alpha_{k_j}\in\{0,1\},\quad 1\leq k_j\leq n\!-\!1.\nonumber
	\end{eqnarray}
\end{small}

where, $d_{\sigma_l,\alpha}$ denotes the distance between the node $x_{\sigma_l}$ and depot $\alpha$. Fig. \ref{fig:mTSP_prob}c shows a schematic of a single-depot returning $m$TSP.

\subsection{Close Enough Traveling Salesmen Problem (CETSP)}\label{subsec:Problem_close}
In this problem, we are given a set of $n$ nodes ($\{x_i\}$), each with a specified radius ($\{\rho_i\}$), and set of $m$ salesmen with the objective of determining the optimal tour such that at least one salesman comes within $\rho_i$ distance from each node $x_i$. CETSPs are used to represent problems such as aerial reconnaissance \cite{mennell2009heuristics} and establishing a wireless meter reader \cite{gulczynski2006close}. The CETSP variant may be applied to any of the TSP class of problems. The most significant difference between point-based TSPs and the CETSP is that due to the radius associated with each node, the CETSP does not define a specific edge between a pair of nodes, rather there is a continuous field of possible edges between a pair of nodes. As a result, there are infinitely many possible solutions to this problem. A single salesman returning CETSP is mathematically described as
\begin{small}
	\begin{eqnarray}\label{eq:opt4}
		\min_{\{v_{ij}\},\{y_j\}}\sum_{j=1}^{n}\left\{\sum_{i=1}^{n}v_{ij}d_{CE}(x_i,y_j)\!+\!d(y_j,y_{j\!+\!1})\right\}; y_{n+1}=y_1 &\nonumber \\
		\text{s.t.}\quad v_{ij}\in\{0,1\},\quad\sum_{i=1}^{n}v_{ij}=1\forall j,\quad\sum_{j=1}^{n}v_{ij}=1,\forall i &\nonumber\\
		\text{where,}\quad d_{CE}(x_i,y_j) = \left\{
				\begin{array}{ll}
					0 & \mbox{if} \|y_j-x_i\|<\rho_i \\
					(\|y_j-x_i\|-\rho_i)^2 & \mbox{else}
				\end{array}
			\right.\nonumber
	\end{eqnarray}
\end{small}
Fig. \ref{fig:mTSP_prob}d shows a schematic of a single salesman returning CETSP.

\section{DETERMINISTIC ANNEALING ALGORITHM:A MAXIMUM ENTROPY PRINCIPLE APPROACH}\label{sec:MEP}
At its core, the deterministic annealing (DA) algorithm solves a facility location problem (FLP): For given $n$ customer locations, find $K$ facility locations such that the total {\em weighted sum of the distance of each customer to its nearest facility is minimized}. In other words, if $x_i$ and $y_j\in\mathbb{R}^N$ denote the locations of $i^{th}$ customer and $j^{th}$ facility, respectively, then the FLP addresses the following optimization problem
\begin{small}
	\begin{equation}\label{eq:FLP}
		\min\limits_{y_j\in\Omega,1\leq j\leq K}\sum\limits_{i=1}^{n}\left\{\min\limits_{y_j,1\leq j\leq K}d(x_i,y_j)\right\},
	\end{equation}
\end{small}
where $d(x_i,y_j)\in\mathbb{R}_+$ denotes the distance between the $i^{th}$ customer location $x_i$ and $j^{th}$ facility location $y_j$, and $\Omega\subset\mathbb{R}^N$ is a compact domain. The solution to an FLP essentially results in a set of clusters $\mathcal{Y} = \{y_j\}$, where a facility $j$ is located at $y_j$ of the $j^{th}$ cluster $C_j$ and each customer is associated only to its nearest facility. Most algorithms for FLP (such as Lloyd's \cite{lloyd1982least}) are overly sensitive to the initial facility locations. This is primarily due to the distributed aspect of FLPs, where any change in the location of the $i^{th}$ customer affects $d(x_i,y_j)$ only with respect to the {\em nearest} facility $j$. The DA algorithm suggested by Rose \cite{rose1998deterministic}, overcomes this sensitivity by allowing each customer to be {\em partially} associated to every cluster through an association probability.

An instance of an FLP is given by the set of facility locations $\mathcal{Y} = \{y_j\}$, and a partition via the set of associates $\mathcal{V} = \{v_{ij}\}$, where
\begin{small}
	\begin{equation}\label{eq:V}
		v_{ij} = \left\{
				\begin{array}{ll}
					1  & \mbox{if } x_i\in C_j \\
					0 & \mbox{else }
				\end{array}
			\right.
	\end{equation}
\end{small}
Borrowing the definition from the data compression literature \cite{cover2012elements}, we associate a {\em distortion} as a measure of the distance of a customer to its nearest facility, with every instance of an FLP, given by
\begin{small}
	\begin{equation}\label{eq:Distortion}
		D(\mathcal{Y},\mathcal{V}) = \sum\limits_{i=1}^{n}\sum\limits_{j=1}^{K}v_{ij}d(x_i,y_j),
	\end{equation}
\end{small}
which is the distortion of the specific, hard-clustering solution. We define the expected distortion as $D = \left\langle D(\mathcal{Y},\mathcal{V})\right\rangle = \sum\limits_{\mathcal{Y},\mathcal{V}}P(\mathcal{Y},\mathcal{V})D(\mathcal{Y},\mathcal{V})$, where $P(\mathcal{Y},\mathcal{V})$ is the instance probability distribution. Since we have no prior knowledge on $P(\mathcal{Y},\mathcal{V})$, we apply the maximum-entropy-principle (MEP) to estimate them. The Shannon entropy term $H = -\sum\limits_{\mathcal{Y},\mathcal{V}}P(\mathcal{Y},\mathcal{V})\log{P(\mathcal{Y},\mathcal{V})}$, widely used in data compression literature \cite{cover2012elements},measures uncertainties in facility locations with respect to the known customer locations. Thus, maximizing the entropy is equivalent to
decreasing the {\em local} influence. The trade-off between minimizing the expected distortion $D$ and maximizing the Shannon entropy $H$ \cite{rose1998deterministic} is achieved by seeking $P(\mathcal{Y},\mathcal{V})$, that minimize the Lagrangian $(\beta D - H)$, where $\beta$ is the Lagrange multiplier (also referred as {\em annealing} parameter). This yields a {\em Gibbs distribution} as
\begin{small}
	\begin{equation}\label{eq:Gibbs}
		P(\mathcal{Y},\mathcal{V}) = \frac{e^{-\beta D(\mathcal{Y},\mathcal{V})}}{\sum\limits_{\mathcal{Y'},\mathcal{V'}}e^{-\beta D(\mathcal{Y'},\mathcal{V'})}}.
	\end{equation}
\end{small} 
Since we are interested in estimating the most probable set of facility locations, we consider the marginal probability, given by
\begin{small}
	\begin{equation}\label{eq:Marginal}
		P(\mathcal{Y}) = \sum\limits_{\mathcal{V}}P(\mathcal{Y},\mathcal{V}) = \frac{Z(\mathcal{Y})}{\sum\limits_{\mathcal{Y'}}Z(\mathcal{Y'})},
	\end{equation}
\end{small}
where $Z(\mathcal{Y}) = \prod\limits_{i=1}^{n}\sum\limits_{j=1}^{K}e^{-\beta d(x_i,y_j)}$. The marginal probability can be rewritten in an explicit Gibbs form as
\begin{small}
	\begin{equation}\label{eq:Marginal_Gibbs}
		P(\mathcal{Y}) = \dfrac{e^{-\beta F(\mathcal{Y})}}{\sum\limits_{\mathcal{Y'}}e^{-\beta F(\mathcal{Y'})}},
	\end{equation}
\end{small}
where $F(\mathcal{Y})$ is the analog of {\em free energy} in statistical mechanics and is given by,
\begin{small}
	\begin{equation}\label{eq:Free_energy}
		F(\mathcal{Y}) = -\frac{1}{\beta}\log{Z(\mathcal{Y})} = -\frac{1}{\beta}\sum\limits_{i=1}^{n}\log\left(\sum\limits_{j=1}^{K}e^{-\beta d(x_i,y_j)}\right).
	\end{equation}
\end{small}
In the DA algorithm, this free energy function is then deterministically optimized at successively increased $\beta$ values over repeated iterations. The set $\mathcal{Y}$ of vectors that optimizes the free energy at each $\beta$ satisfies
\begin{small}
	\begin{equation}\label{eq:Optimal_Y}
		\frac{\partial}{\partial y_j}F = 0 \quad\forall j \quad\Rightarrow \sum\limits_{i=1}^{n}p(j|i)\frac{\partial}{\partial y_j}d(x_i,y_j) = 0 \quad\forall j,
	\end{equation}
\end{small}
where $p(j|i) = \dfrac{e^{-\beta d(x_i,y_j)}}{\sum\limits_{k=1}^{K}e^{-\beta d(x_i,y_k)}}$. The readers are encouraged to refer to \cite{parekh2015deterministic} for detailed analysis on the complexity of the DA algorithm. For implementation on very large datasets, a scalable modification of the DA is proposed in \cite{sharma2006scalable}.

\section{METHODOLOGY: MODIFICATIONS OF THE DA ALGORITHM}\label{sec:methodology}
In this section, we develop a DA based generalized heuristic for variants on the classical TSP. The framework in the DA algorithm is modified to include routing as constrained clustering. Rose has previously explored the application of DA to the TSP \cite{rose1990deterministic}. The heuristic behind the DA based TSP approach is that if we employ same number of facilities as the number of nodes, i.e. $K = n$ and let $\beta\to\infty$, then each node becomes a potential cluster, and the cluster locations $\{y_{j}\}$ coincide with the node locations $\{x_i\}$. The distortion function is
modified to include the tour length. This requires then a second Lagrange multiplier $\theta$ for the tour length component of the distortion function, in addition to the first Lagrange multiplier $\beta$ for the original component of the distortion function. Solving through the gradual change in both Lagrange multipliers leads to a solution of the TSP (and its variants). We also discuss an effective scheme to vary the Lagrange multipliers $\beta$ and $\theta$ in this section.

We now mathematically formulate the problem setting. We are given a set of $n$ nodes whose locations are given by $\{x_i\}\subset\mathbb{R}^2,\quad 1\leq i\leq n$ and a depot with coordinates $\alpha\in\mathbb{R}^2$. These nodes have to be traversed by a maximum of $m$ salesmen under several constraints (which essentially constitute the variants on the TSP) on the optimal tours. As before, we use $\mathcal{Y} = \{y_j\}_{j=1}^n$ and $\mathcal{V} = \{v_{ij}\}_{i,j=1}^n$ to denote the set of facility locations (also referred as {\em codevectors}) and set of associates, respectively. The distance function $d(x_i,y_j)$ between node $i$ and codevector $j$ is considered to be squared-euclidean, i.e., $d(x_i,y_j) = \|x_i-y_j\|_2^2$.

{\bf Remark:} For brevity and ease of exposition, the results are derived for a special case of $m=2$ salesmen. The results are easily extendable for any general $m$ and are stated without any proofs.

\subsection{DA for Non-returning mTSP}\label{subsec:Method_mTSP_nonR}
\begin{figure*}[!t]
	\begin{center}
	\begin{tabular}{ccc}
		\includegraphics[width=0.58\columnwidth]{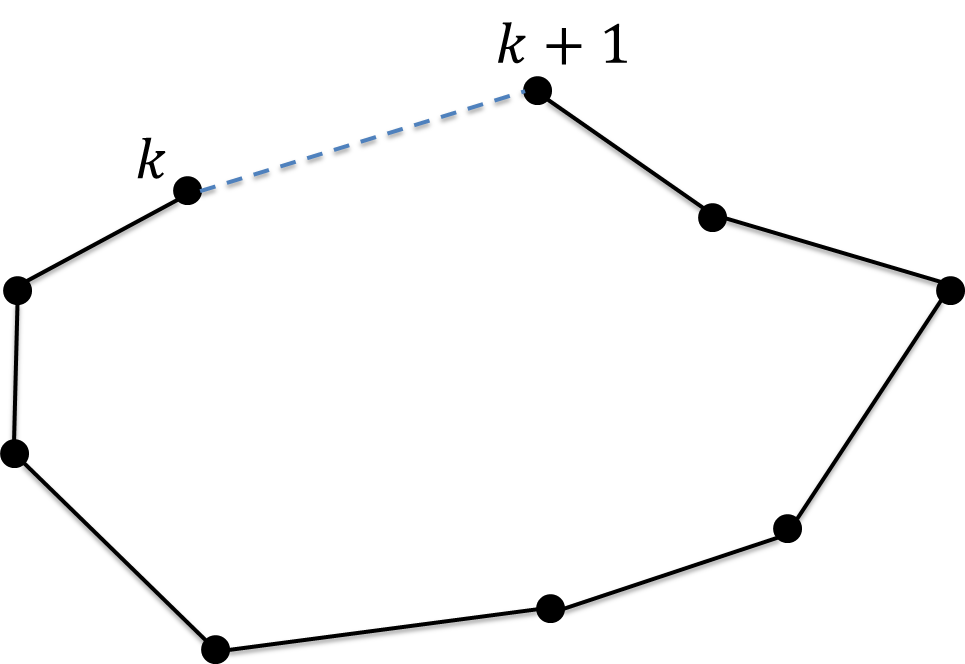}&
		\includegraphics[width=0.58\columnwidth]{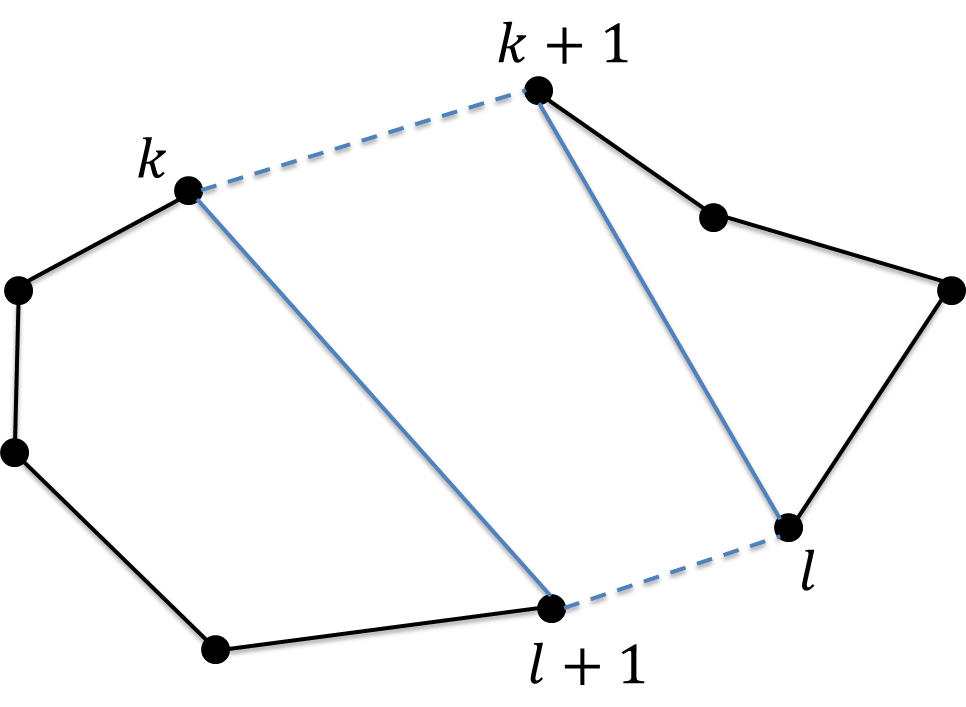}&	
		\includegraphics[width=0.58\columnwidth]{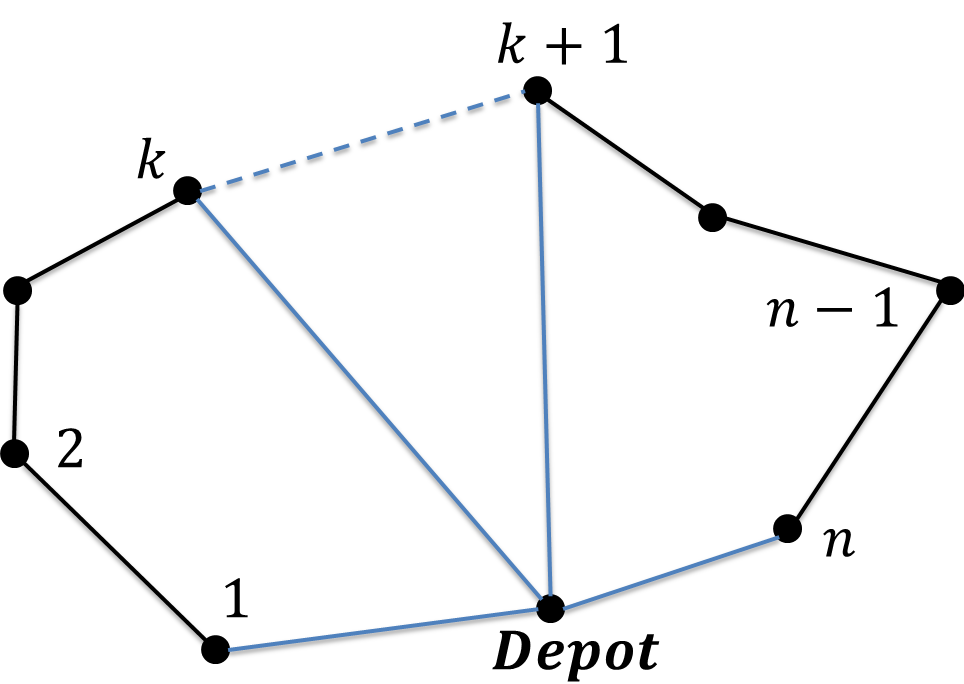}\cr
		(a)&(b)&(c)
	\end{tabular}
	\vspace{-0.5em}
	\caption{{\small Schematic of a (a) Non-returning $2$TSP, with $\mathcal{R}=k$. (b) Returning $2$TSP, with $\mathcal{R}=\{k,l\}$. (c) Returning $2$TSP (with Depot), with $\mathcal{R}=k$. The dashed blue lines indicate the removal of links, while solid blue lines indicate the addition of links.}}
	\label{fig:mTSP_schematic}
	\vspace{-2em}
	\end{center}
\end{figure*}
In this case, an instance is defined by three parameters, $\mathcal{Y}, \mathcal{V}$ and $\mathcal{R}$. Here $\mathcal{Y}$ corresponds to the set of codevectors $y_j$, $\mathcal{V}=\{v_{ij}\}$ is a set of associates and describes the membership of node $x_i$ to codevector $y_j$ (see Eq. \ref{eq:V}), and $\mathcal{R}$ is a set of locations of the partition representing the breaks between subsequent salesmen in the chain of consecutive codevectors (see Fig. \ref{fig:mTSP_schematic}a). We first consider the case for $m = 2$ salesmen. Therefore an instance $\mathcal{R}$ is mathematically described as
\begin{small}
	\begin{equation*}
		\mathcal{R} = k; \quad \mbox{if there is no link b/w $y_k$ and $y_{k+1}$}.
	\end{equation*}
\end{small}
As this is the non-returning version of the mTSP, there is no connection between $y_1$ and $y_n$. Therefor $y_0 = y_{n+1} = 0$ when they do appear in the equations. For a given instance of the problem $(\mathcal{Y},\mathcal{V},\mathcal{R})$, the distortion function in (\ref{eq:Distortion}) is modified as
\begin{small}
	\begin{equation}\label{eq:D1_nonR_2TSP}
		D(\mathcal{Y},\mathcal{V},\mathcal{R}) = D_1(\mathcal{Y},\mathcal{V}) + D_2(\mathcal{Y}) + D_3(\mathcal{Y},\mathcal{R}),
	\end{equation}
\end{small}
where $D_1(\mathcal{Y},\mathcal{V})$ is same as defined in (\ref{eq:Distortion}). $D_2(\mathcal{Y})$ captures the tour length in the cost function to represent a TSP and is given by
\begin{small}
	\begin{equation}\label{eq:D2_nonR_2TSP}
		D_2(\mathcal{Y}) = \theta\sum\limits_{j=1}^{n-1}d(y_j,y_{j+1}).
	\end{equation}
\end{small}
The final component $D_3(\mathcal{Y},\mathcal{R})$ represents the partition of codevectors for the independent salesmen and subtracts the distance at the partition between the codevector $y_k$ and $y_{k+1}$ from the original distortion function.
\begin{small}
	\begin{equation}\label{eq:D3_nonR_2TSP}
		D_3(\mathcal{Y},\mathcal{R}) = -\theta d(y_k,y_{k+1})
	\end{equation}
\end{small}
Similar to the original formulation of the DA algorithm, the probability of any instance $(\mathcal{Y},\mathcal{V},\mathcal{R})$ is determined by the MEP. Following the steps adopted in Sec. \ref{sec:MEP}, the free energy of this system is obtained as
\begin{small}
	\begin{eqnarray}\label{eq:F_nonR_2TSP}
		F &=& -\frac{1}{\beta}\sum\limits_{i=1}^{n}\log\left(\sum\limits_{j=1}^{n}e^{-\beta d(x_i,y_j)}\right) + \theta\sum\limits_{j=1}^{n-1}d(y_j,y_{j+1})\nonumber \\
		&&-\frac{1}{\beta}\log\left(\sum\limits_{k=1}^{n-1}e^{\beta\theta d(y_k,y_{k+1})}\right).
	\end{eqnarray}
\end{small}
Taking the derivative of (\ref{eq:F_nonR_2TSP}) with respect to each codevector allows determination of the set of codevectors that maximize entropy in the system.
\begin{small}
	\begin{eqnarray}\label{eq:dF_nonR_2TSP}
		&&\frac{\partial F}{\partial y_j} = -2\sum\limits_{i=1}^{n}p(j|i)(y_j-x_i) + 2\theta(2y_j-y_{j+1}-y_{j-1}) \nonumber \\
		&&+ 2\theta \left(y_{j+1}-y_j\right)\mathcal{P}(j) +  2\theta\left(y_{j-1}-y_j\right)\mathcal{P}(j-1) = 0\nonumber
	\end{eqnarray}
\end{small}
where $p(j|i)$ is same as before (see Eq. (\ref{eq:Optimal_Y})). $\mathcal{P}(j)$ represents the probability that the partition occurs at codevector $j$ (i.e. $\mathcal{R}=j$) and is given by
\begin{small}
	\begin{equation}\label{eq:pr_nonR_2TSP}
		\mathcal{P}(j) = \dfrac{e^{\beta\theta d(y_j,y_{j+1})}}{\sum\limits_{j=1}^{n-1}e^{\beta\theta d(y_j,y_{j+1})}}.
	\end{equation}
\end{small}
Solving for each $y_j$ provides the solution to the system at this pair of $\beta$ and $\theta$ values, so that for every codevector
\begin{small}
	\begin{equation}\label{eq:Opt_Y_nonR_2TSP}
		y_j = \frac{\sum\limits_{i=1}^{n}p(j|i)x_i + \theta y_{j+1}(1\!-\!\mathcal{P}(j)) + \theta y_{j-1}(1\!-\!\mathcal{P}(j\!-\!1))}{\sum\limits_{i=1}^{n}p(j|i) + \theta\left(2\!-\mathcal{P}(j)\!-\mathcal{P}(j\!-\!1)\!\right)}
	\end{equation}
\end{small}
Note that the Eq. \ref{eq:Opt_Y_nonR_2TSP} is only slightly more complex than the basic TSP proposed in \cite{rose1990deterministic}. In fact, setting $\mathcal{P}(j) = 0,\forall j$ transforms the $2$TSP into the basic TSP formulation.

\underline{\em Extension to $m$TSP for a general $m\geq 2$}: In this case, the partition set $\mathcal{R}$ contains $m\!-\!1$ points, i.e., $\mathcal{R}\!=\!\{k_1,k_2,...,k_{m\!-\!1}\}$. The distortion function $D_3(\mathcal{Y},\mathcal{R})$ in (\ref{eq:D3_nonR_2TSP}) and the probability of a partition $\mathcal{P}(k_1,..,k_{m\!-\!1})$ in (\ref{eq:pr_nonR_2TSP}) are modified as
\begin{small}
	\begin{eqnarray}\label{eq:D3_Pr_Y_nonR_mTSP}
		D_3(\mathcal{Y},\mathcal{R}) &=& -\theta\sum_{i=1}^{m-1}d(y_{k_i},y_{k_i\!+\!1})\nonumber\\
		\mathcal{P}(k_1,..,k_{m\!-\!1}) &=& \frac{e^{\beta\theta\sum_{i=1}^{m-1}d(y_{k_i},y_{k_i\!+\!1})}}{\sum\limits_{\mathcal{R}}e^{\beta\theta\sum_{i=1}^{m-1}d(y_{k_i},y_{k_i\!+\!1})}}.\nonumber
	\end{eqnarray}
\end{small}
Note that the probability $\mathcal{P}(k_1,..,k_{m\!-\!1})$ is symmetric in its arguments, i.e.,\\
\centerline{$\mathcal{P}(k_1,k_2,..,k_{m\!-\!1}) = \dots = \mathcal{P}(k_{m\!-\!1},k_1,..,k_{m\!-\!2})$}.\\
Let us define $\mathcal{\tilde{R}} = \mathcal{R}\setminus\{k_1\}$ and denote $\mathcal{P}(j,k_2,..,k_m)$ by $\mathcal{P}(j,\mathcal{\tilde{R}})$. Then the characteristic equation for each codevector is given by
\begin{small}
	\begin{equation*}\label{eq:Opt_Y_nonR_mTSP}
		y_j = \frac{\splitfrac{\sum\limits_{i=1}^{n}p(j|i)x_i + \theta y_{j+1}\bigg(1\!-\!(m\!-1)\!\sum\limits_{\mathcal{\tilde{R}}}\mathcal{P}(j,\mathcal{\tilde{R}})\bigg)}{+\theta y_{j-1}\bigg(1\!-\!(m\!-1)\!\sum\limits_{\mathcal{\tilde{R}}}\mathcal{P}(j\!-\!1,\mathcal{\tilde{R}})\bigg)}}{\sum\limits_{i=1}^{n}p(j|i)+\theta\Big(2\!-\!(m\!-\!1)\sum\limits_{\mathcal{\tilde{R}}}\Big(\mathcal{P}(j,\mathcal{\tilde{R}})+\mathcal{P}(j-1,\mathcal{\tilde{R}})\Big)\Big)}
	\end{equation*}
\end{small}

\subsection{DA for Returning mTSP}\label{subsec:Method_mTSP_R}
In case of returning $m$TSP, the start and end positions of each salesmen must be coincident. In this case, the partition function $D_3(\mathcal{Y},\mathcal{R})$ not only considers the distance between the codevectors where the partition occurs, it must also account for the distance incurred in completing the continuous tour by reconnecting to the other end of the loop (see Fig. \ref{fig:mTSP_schematic}b). Similar to the non-returning $m$TSP, we first derive the results for $m=2$ salesmen and then later extend it for a general $m$. The partition parameter $\mathcal{R}$ in this case is described by two parameters.
\begin{small}
	\begin{equation*}
		\mathcal{R} = k,l\quad \mbox{if}\bigg\{ \begin{array}{c}
			\text{no links b/w $y_k$ and $y_{k+1}$, $y_l$ and $y_{l+1}$};\\
			\text{links b/w $y_k$ and $y_{l+1}$, $y_l$ and $y_{k+1}$};
		\end{array}
	\end{equation*}
\end{small}
It should be noted that the codevectors $y_1$ and $y_n$ are considered to be adjacent, i.e. $y_0 = y_n$ and $y_{n+1} = y_1$. The tour-length distortion function is given by $D_2(\mathcal{Y})=\sum\limits_{j=1}^{n}d(y_j,y_{j+1})$. The distortion function $D_3(\mathcal{Y},\mathcal{R})$ pertaining to the partition parameter is defined as
\begin{small}
	\begin{equation*}\label{eq:D3_R_2TSP}
		D_3 = \theta\big(\!-\!d(y_k,y_{k+1})\!-\!d(y_l,y_{l+1})\!+\!d(y_k,y_{l+1})\!+\!d(y_l,y_{k+1})\!\big).
	\end{equation*}
\end{small}
An important consequence of this framework is that if $k=l$, then the problem reduces to the classical returning TSP. Thus, this framework allows automatic determination of the optimal number of salesmen. With the above modified distortion functions, the update equation for each codevector $y_j$ is given by
\begin{small}
	\begin{equation}\label{eq:Opt_Y_R_2TSP}
		y_j = \frac{\splitfrac{\sum\limits_{i=1}^{n}p(j|i)x_i+2\theta\Big(\sum\limits_{l\neq j-1}\mathcal{P}(j,l)y_{l+1}+\sum\limits_{l\neq j}\mathcal{P}(j\!-\!1,l)y_{l}\Big)}{+\theta\left[(1\!-\!2\sum_{l}\mathcal{P}(j,l))y_{j+1}\!+\!(1\!-\!2\sum_{l}\mathcal{P}(j\!-\!1,l))y_{j-1}\right]}}{\sum_{i=1}^{n}p(j|i)+2\theta\big(1-2\mathcal{P}(j,j-1)\big)}
	\end{equation}
\end{small}
\underline{\em Extension to $m$TSP for a general $m\geq 2$}: We now extend out formulation for returning $m$TSP to a general $m\geq 2$. The partition set $\mathcal{R}$ is defined as
\begin{small}
	\begin{equation*}
		\mathcal{R} = \{k_1,k_2\dots k_m\} \bigg\{ \begin{array}{c}
			\text{no links b/w $y_{k_i}$ and $y_{k_i+1}$};\\
			\text{links b/w $y_{k_i}$ and $y_{k_{(i-1)}+1}$};
		\end{array}
	\end{equation*}
\end{small}
with the understanding that $k_0 = k_m$. The distortion function and the corresponding probability distribution pertaining to the partition parameter are given by
\begin{small}
	\begin{eqnarray}\label{eq:D3_P_R_mTSP}
		D_3(\mathcal{Y},\mathcal{R}) &=& \theta\sum_{\mathcal{R}}\left\{d(y_{k_i},y_{k_{(i-1)}+1})-d(y_{k_i},y_{k_i+1})\right\}\nonumber\\
		\mathcal{P}(\mathcal{R}) &=& \frac{\sum_{i=1}^{m}e^{-\beta\theta\left\{d(y_{k_i},y_{k_{(i-1)}+1})-d(y_{k_i},y_{k_i+1})\right\}}}{\sum_{\mathcal{R}}\sum_{i=1}^{m}e^{-\beta\theta\left\{d(y_{k_i},y_{k_{(i-1)}+1})-d(y_{k_i},y_{k_i+1})\right\}}}. \nonumber
	\end{eqnarray}
\end{small}
Please note that the probability distribution $\mathcal{P}(k_1,\dots ,k_m)$ is symmetric in its arguments. Let us define $\mathcal{\tilde{R}}=\mathcal{R}\setminus\{k_1,k_2\}$. Then the corresponding update equation for each codevector $y_j$ is given by
\begin{small}
	\begin{equation*}\label{eq:Opt_Y_R_mTSP}
		y_j = \frac{\splitfrac{m\theta\sum_{\tilde{R}}\mathcal{P}(j,j\!-\!1,\tilde{R})(y_{j\!+\!1}\!+\!y_{j\!-\!1})\!+\!\sum_{i=1}^{n}p(j|i)x_i} {\splitfrac{+\!m\theta\sum_{k\neq m\!-\!1}^{n}\sum_{\mathcal{\tilde{R}}}\mathcal{P}(j,k,\mathcal{\tilde{R}})(y_{k\!+\!1}-y_{j\!+\!1}) } {+\!m\theta\sum_{k\neq m}^{n}\sum_{\mathcal{\tilde{R}}}\mathcal{P}(j\!-\!1,k,\mathcal{\tilde{R}})(y_{k}-y_{j\!-\!1})}} } {\sum_{i\!=\!1}^{n}p(j|i)\!+\!\theta\left\{1\!-\!2m\left[\sum_\mathcal{\tilde{R}}\mathcal{P}(j,j\!-\!1,\mathcal{\tilde{R}})\right]\right\}}.
	\end{equation*}
\end{small}

\subsection{DA for Single Depot Returning mTSP}\label{subsec:Method_mTSP_RD}
Fig. \ref{fig:mTSP_schematic}c shows the schematic of the proposed framework for the returning version of a single depot multiple salesmen problem. We denote the depot location by $\alpha\in\mathbb{R}^2$. We first formulate the framework for $m=2$ salesmen. The partition parameter $\mathcal{R}$ in this case is defined as
\begin{small}
	\begin{equation*}
		\mathcal{R} = k, \quad \Big\{\begin{array}{c}
			\text{no link b/w $y_k$ and $y_{k\!+\!1}$};\\
			\text{links b/w $y_k$ and $\alpha$, $y_{k\!+\!1}$ and $\alpha$}.
		\end{array}
	\end{equation*}
\end{small}
The distortion function $D_3(\mathcal{Y},\mathcal{R})$ and the corresponding probability distribution $\mathcal{P}(\mathcal{R}=k)$ pertaining to the partition parameter are given by
\begin{small}
	\begin{eqnarray}
		D_3(\mathcal{Y},\mathcal{R}) &=& \theta\left(-d(y_k,y_{k\!+\!1})\!+\!d(y_k,\alpha)\!+\!d(y_{k\!+\!1},\alpha)\right)\nonumber\\
		P(\mathcal{Y}) &=& \frac{ e^{-\beta\theta\{-d(y_k,y_{k\!+\!1})\!+\!d(y_k,\alpha)\!+\!d(y_{k\!+\!1},\alpha)\}} }{\sum_{k=1}^{n-1}e^{-\beta\theta\{-d(y_k,y_{k\!+\!1})\!+\!d(y_k,\alpha)\!+\!d(y_{k\!+\!1},\alpha)\}}}.\nonumber
	\end{eqnarray}
\end{small}
The distortion function corresponding to the tour-length constraint is modified to include the links between $y_1$ and $\alpha$, and between $y_n$ and $\alpha$, i.e., $D_2(\mathcal{Y})=\theta\left\{\sum_{j=1}^{n\!-\!1}d(y_j,y_{j\!+\!1})\!+\!d(y_1,\alpha)\!+\!d(y_n,\alpha)\right\}$. If we define $\mathcal{P}(0)=\mathcal{P}(n)=1$, then the corresponding update equation for each codevector $y_j$ is given by
\begin{small}
	\begin{equation*}
		y_j = \frac{\splitfrac{\sum_{i=1}^{n}p(j|i)x_i + \theta\left\{\mathcal{P}(j)\!+\!\mathcal{P}(j\!-\!1)\right\}\alpha} {+\theta\left\{1-\mathcal{P}(j)\right\}y_{j\!+\!1} + +\theta\left\{1-\mathcal{P}(j\!-\!1)\right\}y_{j\!-\!1}} } {2\theta+\sum_{i=1}^{n}p(j|i)}.
	\end{equation*}
\end{small}
{\bf Remark}: While the proposed framework finds optimal location where the partition occurs, the resulting optimal tour-lengths for the two salesmen can differ considerably. Due to various physical constraints such as fuel and capacity constraints in single-depot vehicle routing problems, it is often required that the tour-lengths for the two salesmen must be comparable, without significantly affecting the total distance traveled by the two salesmen. This constraint can be incorporated into the current framework by appropriately modifying the distortion functions as
\begin{small}
	\begin{eqnarray}
		D'_2(\mathcal{Y}) &=& D_2+\eta\theta(d(y_1,\alpha)-d(y_n,\alpha))\nonumber\\
		D'_3(\mathcal{Y},\mathcal{R}) &=& D_3\!+\!\eta\theta\left\{\sum_{i=1}^{k}d(y_i,y_{i\!+\!1})-\sum_{i=k\!+\!1}^{n}d(y_i,y_{i\!+\!1})\right\}\nonumber
	\end{eqnarray}
\end{small}
where $\eta$ is a constant parameter and captures the trade-off between between obtaining the optimal total distance and minimizing the imbalance between the two tour-lengths. The probability distribution and the update equation for codevectors are appropriately modified.

\underline{\em Extension to $m$TSP for general $m\geq 2$}: In this case, the partition set $\mathcal{R}$ consists of $m\!-\!1$ points, i.e., $\mathcal{R}=\{k_1,\dots ,k_{m\!-\!1}\}$. The distortion and the probability distribution pertaining to the partition parameter are modified as
\begin{small}
	\begin{eqnarray}
		D_3(\mathcal{Y},\mathcal{R}) = \theta\sum_{i=1}^{m-1}\left\{-d(y_{k_i},y_{k_i\!+\!1})\!+\!d(y_{k_i},\alpha)\!+\!d(y_{k_i\!+\!1},\alpha)\right\}\nonumber\\
		\mathcal{P} = \frac{e^{-\beta\theta\sum_{i=1}^{m-1}\left\{-d(y_{k_i},y_{k_i\!+\!1})\!+\!d(y_{k_i},\alpha)\!+\!d(y_{k_i\!+\!1},\alpha)\right\}}} {\sum_\mathcal{R}e^{-\beta\theta\sum_{i=1}^{m-1}\left\{-d(y_{k_i},y_{k_i\!+\!1})\!+\!d(y_{k_i},\alpha)\!+\!d(y_{k_i\!+\!1},\alpha)\right\}}}\nonumber
	\end{eqnarray}
\end{small}
Note that the probability distribution $\mathcal{P}(k_1,\dots ,k_{m\!-\!1})$ is symmetric in its arguments. If we define $\mathcal{\tilde{R}}$ by $\mathcal{\tilde{R}}=\mathcal{R}\setminus k_1$, then the corresponding update equation for each codevector $y_j$ is given by
\begin{small}
	\begin{equation*}
		y_j = \frac{\splitfrac{\sum_{i=1}^{n}p(j|i)x_i+(m\!-\!1)\theta\sum_\mathcal{\tilde{R}}\big(\mathcal{P}(j,\mathcal{\tilde{R}})\!+\!\mathcal{P}(j\!-\!1,\mathcal{\tilde{R}})\big)\alpha} {\splitfrac{+\!\theta\big\{1\!-\!(m\!-\!1)\sum_\mathcal{\tilde{R}}\mathcal{P}(j,\mathcal{\tilde{R}})\big\}y_{j\!+\!1}\!} {+\!\theta\big\{1\!-\!(m\!-\!1)\sum_\mathcal{\tilde{R}}\mathcal{P}(j-1,\mathcal{\tilde{R}})\big\}y_{j\!-\!1}\!}} } {\sum_{i=1}^{n}p(j|i)+2\theta}.
	\end{equation*}
\end{small}
 
\subsection{DA for Close Enough TSP}\label{subsec:Method_close}
In the close enough traveling salesman problem (CETSP), an additional radius parameter ($\rho_i$) corresponding to each node $x_i$ is included in the optimization framework. For ease of exposition, we consider a single salesman CETSP in this work, however, the framework can be modified to additionally incorporate any of the aforementioned variants. Note that there are no partition parameters for a single salesman returning TSP. The distance between the node and the codevector pairs is modified as
\begin{small}
	\begin{equation}\label{eq:D_CE}
		d_{CE}(x_i,y_j,\rho_i) = \left(\|y_j-x_i\|-\rho_i\right)^2
	\end{equation}
\end{small}

The distortion functions corresponding to the node-codevector distances and the tour-length constraints are respectively given by
\begin{small}
	\begin{eqnarray}\label{eq:D1_D2_CETSP}
		D_1(\mathcal{Y},\mathcal{V}) &=& \sum_{i=1}^{n}\sum_{j=1}^{n}v_{ij}d_{CE}(x_i,y_j,\rho_i)\nonumber\\
		D_2(\mathcal{Y}) &=& \theta\sum_{j=1}^nd(y_j,y_{j\!+\!1})
	\end{eqnarray}
\end{small}
The free energy of this system is obtained as
\begin{small}
	\begin{equation*}
		F = -\frac{1}{\beta}\sum_{i=1}^{n}\log\Big(\sum_{j=1}^{n}e^{-\beta d_{CE}(x_i,y_j,\rho_i)}\Big)\!+\!\theta\sum_{j=1}^{n}d(y_j,y_{j\!+\!1}).
	\end{equation*}
\end{small}
Taking derivative of the free-energy term and setting it to $0$, we obtain the update equation for each codevector given by
\begin{small}
	\begin{equation*}
		y_j=\frac{\sum_{i=1}^{n}p(j|i)(x_i+\rho_i\sign{(y_j-x_i)})\!+\!\theta(y_{j\!+\!1}+y_{j\!-\!1})} {2\theta+\sum_{i=1}^{n}p(j|i)},
	\end{equation*}
\end{small}
where, the association probability distribution is now given by $p(j|i) = \Big(\frac{e^{-\beta d_{CE}(x_i,y_j,\rho_i)}}{\sum_{k=1}^{n}e^{-\beta d_{CE}(x_i,y_k,\rho_i)}}\Big)$ and $\sign(\cdot)$ is a vector-valued {\em signum} function.

\subsection{Controlling Lagrange Multipliers}\label{subsec:Lagrange_Control}
It is desirable to have a consistent and repeatable method for varying the Lagrange multipliers $\beta$ and $\theta$ which govern the distortion function and the tour length. In this study, the $\beta$ multiplier is considered as the main driver, and the $\theta$ multiplier is secondary. As such, the $\theta$ parameter is decreased according to an exponential function until a stable tour length is reached, at which point $\beta$ is increased according to an exponential function. This process is repeated until a sufficiently high $\beta$ value and sufficiently low $\theta$ value are both reached, leaving the final solution. This is addressed by Rose \cite{rose1990deterministic} in the context of classical TSP, which is generalized to $m$TSP case in this work.

\section{RESULTS AND DISCUSSIONS}\label{sec:results}
This section provides an overview of the results of the MATLAB implementations (on an Intel $i5-4200$U @ $2.30$GHz machine) of the proposed heuristic. As yet, the MATLAB code used for this implementation has not been optimized for minimum computation time, so valid comparisons on the basis on run time are not currently available, however the heuristic is shown to achieve high quality results based on tour lengths in fairly reasonable amount of time. Tests are performed on both the $m$TSP heuristic and the CETSP variation. The heuristics are evaluated on synthetic data.
\begin{figure*}
	\begin{center}
	\begin{tabular}{ccc}
		\includegraphics[width=0.78\columnwidth]{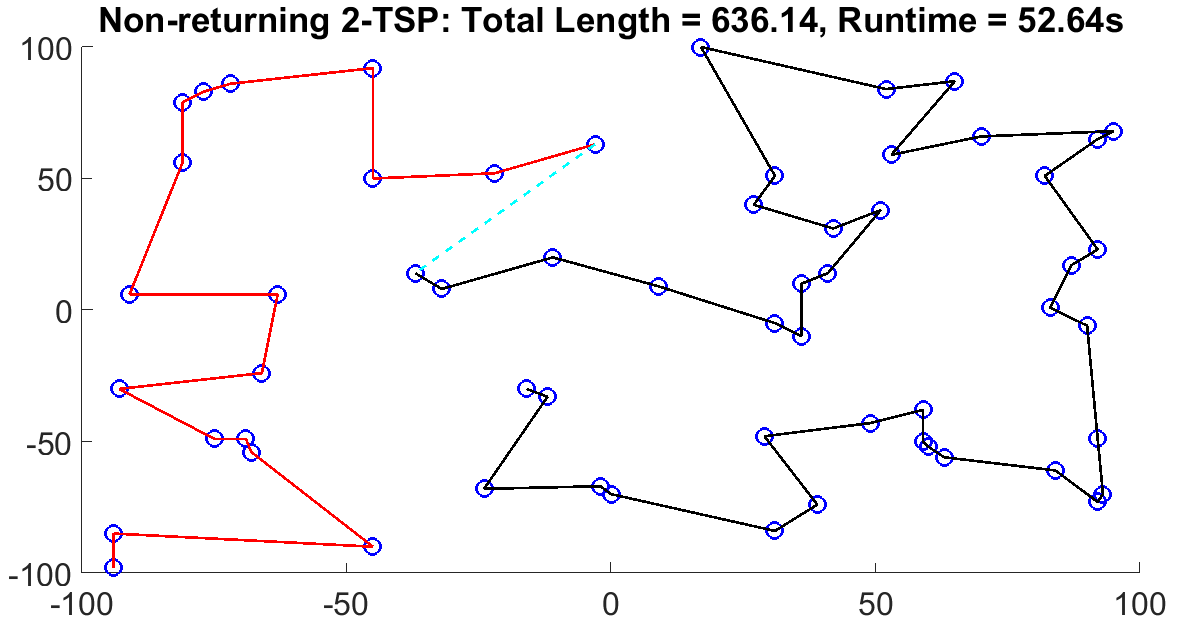}&
		\includegraphics[width=0.75\columnwidth]{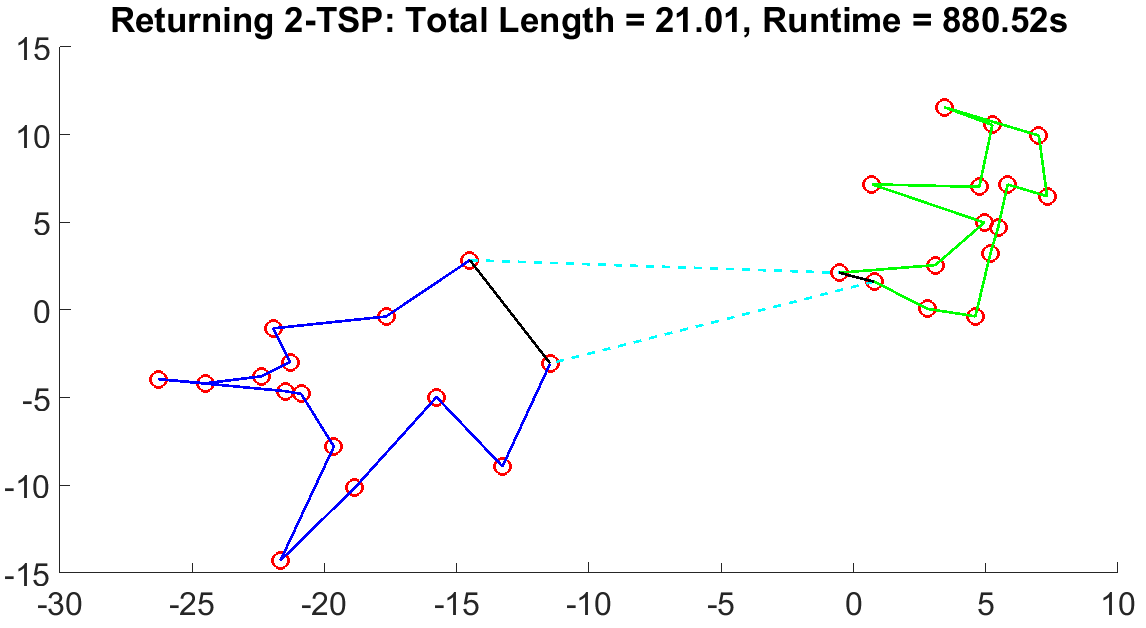}&	
		\includegraphics[width=0.4\columnwidth]{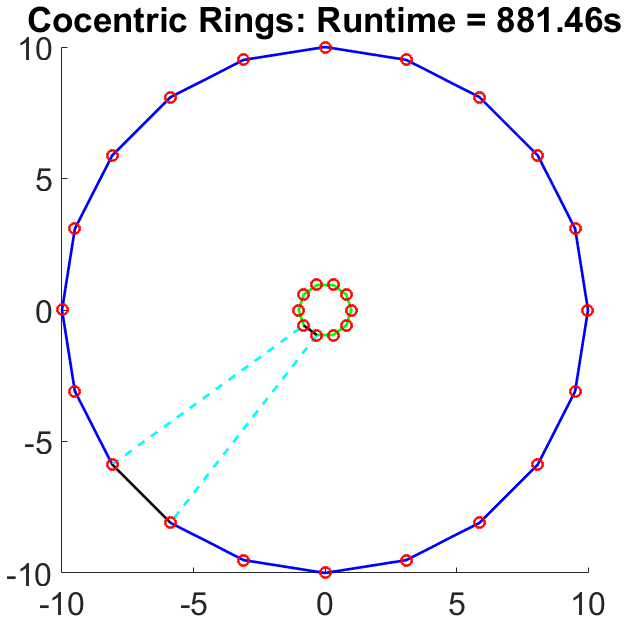}\cr
		(a)&(b)&(c)
	\end{tabular}
	\caption{{\small (a) Result for non-returning $2$TSP for $59$ nodes data. The dashed line indicates the removal of the link. (b) Result for returning $2$TSP for $59$ nodes data. The dashed lines indicate the removal of links, while the black solid lines indicate the addition of links. (c) Returning $2$TSP version for cocentric rings.}}
	\label{fig:mTSP_R_nonR}
	\end{center}
\end{figure*}

\underline{\em Non-returning $2$TSP}: Fig. \ref{fig:mTSP_R_nonR}a shows the non-returning $2$TSP result for a synthetic $59$ nodes data. The tours of the two salesmen are shown in {\em red} and {\em black} respectively, with {\em cyan} dashed line indicating the partition.

\underline{\em Returning $2$TSP}: Fig. \ref{fig:mTSP_R_nonR}b shows the returning $2$TSP result for a randomly generated $30$ nodes data. The heuristic is able to find the two largest links to be removed from the sequence of codevectors. Fig. \ref{fig:mTSP_R_nonR}c shows the returning $2$TSP result for a $30$ nodes cocentric rings arrangement. The DA based heuristic finds the two most optimal routes for this configuration. Note that this dataset is particularly challenging for heuristics such as {\em cluster-first route-second}, where clustering the data first will either result in two symmetric subsets or the only cluster identified will be at the origin and when the two salesmen are allocated to the nodes, there is no way to effectively partition the set into two distinct subsets based on the information provided by the clustering solution.

\begin{figure*}
	\begin{center}
	\begin{tabular}{ccc}
		\includegraphics[width=0.95\columnwidth]{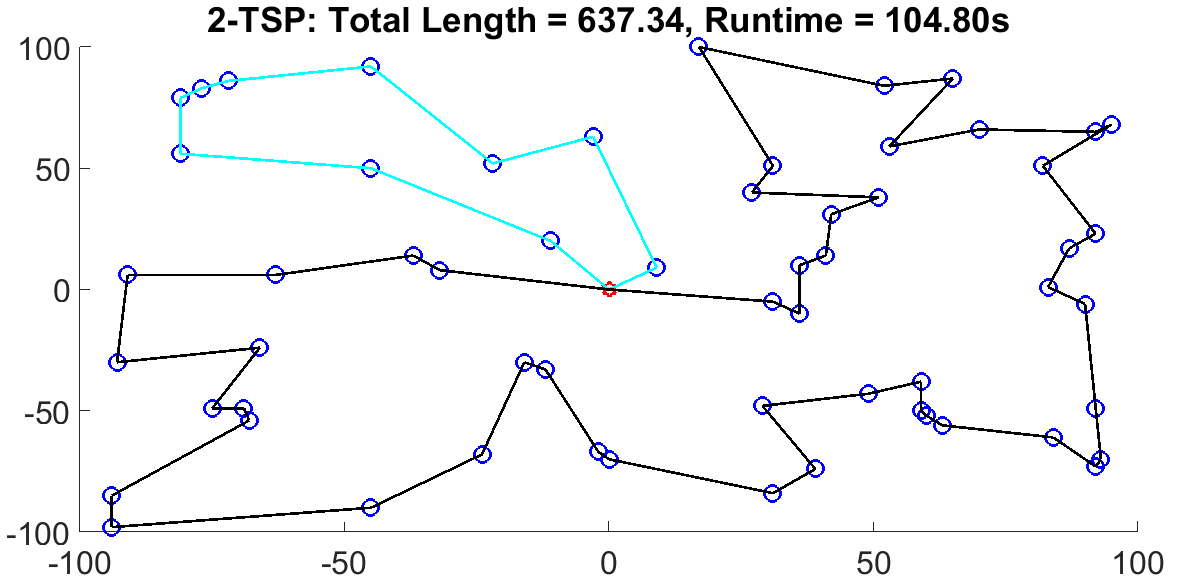}&
		\includegraphics[width=0.9\columnwidth]{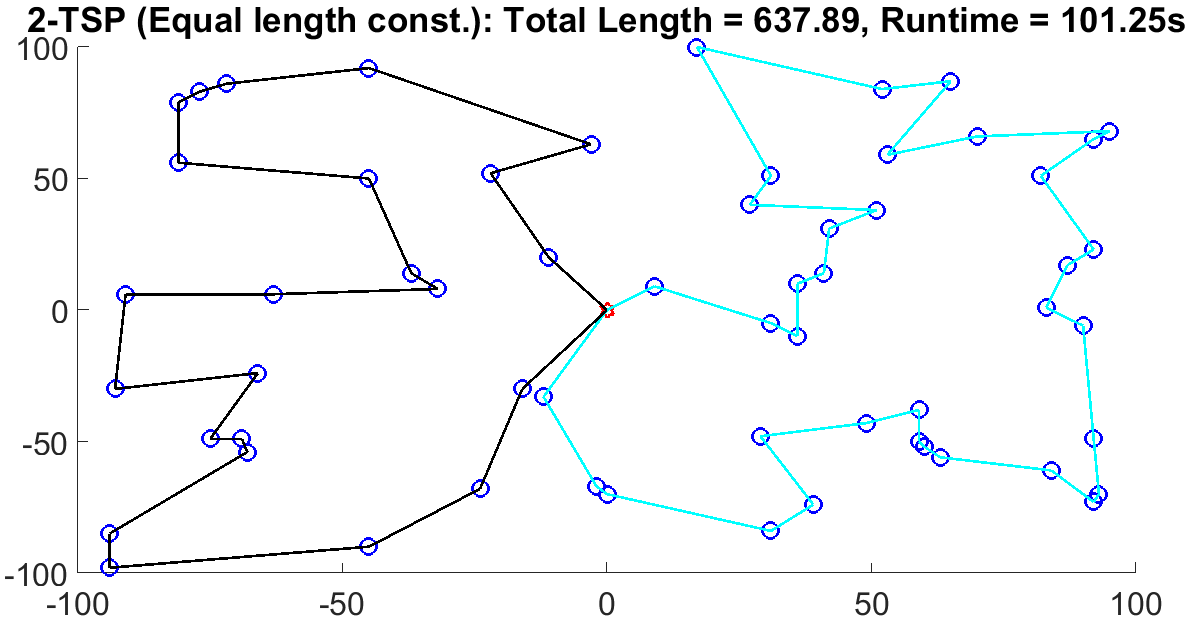}\cr
		(a)&(b)
	\end{tabular}
	\caption{{\small (a) Single-depot returning $2$TSP solution to a $59$ nodes data. The depot location is denoted by {\em red} marker. (b) Solution to the same $59$ nodes data with equality constraint. Note that the optimal length is only slightly bigger than the earlier case.}}
	\label{fig:mTSP_depot}
	\end{center}
\end{figure*}

\underline{\em Single depot returning $2$TSP}: Fig. \ref{fig:mTSP_depot}a shows the implementation results for the $59$ nodes (and a depot) data. The two tours are shown in {\em cyan} and {\em black} colors respectively. Fig. \ref{fig:mTSP_depot}b shows the result with an additional equality constraint. Note that the optimal length in this case is only marginally greater than the optimal length in the previous case.

\begin{figure}
	\begin{center}
		\includegraphics[width=0.95\columnwidth]{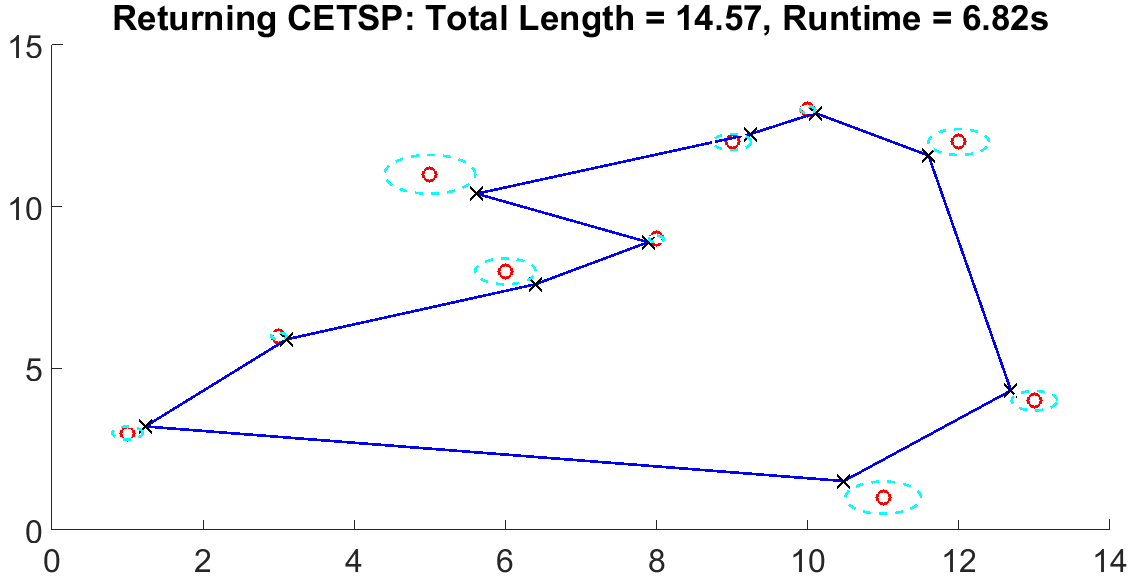}
		\caption{{\small CETSP result for single salesman $10$ nodes returning TSP. The {\em red} markers denote the node locations, while the {\em black} `$\times$' denote the codevector locations. The {\em cyan} circles correspond to the radii $\rho_i$.}}
		\label{fig:CETSP}
	\end{center}
\end{figure}

\underline{\em Returning CETSP}: Fig. \ref{fig:CETSP} shows the implementation results for the CETSP on a randomly generated $10$ nodes data with the additional radius parameter. It is difficult to determine whether the algorithm arrives at an optimal solution because this is much more difficult to check manually and unlike the standard TSP, there is no database of optimal tours for the CETSP. We have compared the heuristic against one of the 100 node sets (kroD100 from TSPLIB \cite{reinelt1991tsplib}) tested by Mennell for equal radii of 11.697\cite{mennell2009heuristics}. Mennell achieves a tour length of 58.54 units with a 0.3 overlap ratio on the data. However, there are no details on the calculation time. The MEP based heuristic finds an optimal tour length of 64.99 units in 949 seconds. Note that in the current formulation, there is a penalty for a codevector existing either inside or outside of the circle. However, according to the problem formulation, there should be no penalty when the codevector exists within the radius of the node. This can be addressed by setting the derivative of the distance function $d_{CE}(x_i,y_j,\rho_i)$ with respect to $y_j$ to zero whenever $y_{j}$ exists within $\rho_i$ distance from the node $x_i$. This negates the penalty incurred for placing a codevector within the radius of a node and should help this heuristic identify more accurate solutions.

\section{CONCLUSIONS AND FUTURE WORKS}\label{sec:conclusion}
In this paper we explore the Maximum-Entropy-Principle as a heuristic for the TSP, as well as many variants. Because the algorithm is independent of the edges between nodes, it has more flexibility to address variants such as the CETSP and the mTSP. The algorithm produces high-quality solutions for some challenging scenarios, such as, cocentric rings.

The next steps for this heuristic framework should be developing the formulation for further variants on the basic TSP. There remain significant opportunities to optimize the code implementation of the MEP framework to achieve more favorable computation times, at which point this algorithm can be run on benchmark mTSP cases and compared against many of the conventional heuristics. Another scope of this work lies in hybridization. It is suggested that the strengths of conventional heuristics should be explored to determine whether any of them complement this heuristic such that a hybrid model could be more successful than either of the individual heuristics.

\section{ACKNOWLEDGMENT}\label{sec:acknowledgement}

The authors would like to acknowledge NSF grants ECCS 15-09302, CMMI 14-63239 and CNS 15-44635 for supporting this work.


\bibliographystyle{IEEEtran} 
\bibliography{myRef}

\end{document}